\documentclass[leqno]{amsart}
\usepackage{amssymb}

\newcommand{\Q}{{\mathbb{Q}}}
\newcommand{\F}{{\mathbb{F}}}
\newcommand{\C}{{\mathbb{C}}}
\newcommand{\Z}{{\mathbb{Z}}}
\newcommand{\R}{{\mathbb{R}}}

\newcommand{\fB}{{\mathfrak{B}}}

\newcommand\Irr{\operatorname{Irr}}
\newcommand\Ind{\operatorname{Ind}}

\newcommand\Res{\operatorname{Res}}

\newtheorem{thm}{Theorem}[section]
\newtheorem{cor}[thm]{Corollary}
\newtheorem{prop}[thm]{Proposition}
\newtheorem{lem}[thm]{Lemma}

\theoremstyle{definition}

\newtheorem{defn}[thm]{Definition}

\newtheorem{rem}[thm]{Remark}

\renewcommand{\leq}{\leqslant}
\renewcommand{\geq}{\geqslant}
\renewcommand{\atop}[2]{\genfrac{}{}{0pt}{}{#1}{#2}}

\address{Institut Girard Desargues, Universit\'e Lyon 1, 21 av Claude 
Bernard, 69622 Villeurbanne cedex, France}

\email{geck@desargues.univ-lyon1.fr}

\begin{document}

\title{The Schur indices of the cuspidal unipotent characters of the
finite Chevalley groups $E_7(q)$}

\author{Meinolf Geck}

\subjclass[2000]{Primary 20C15; Secondary 20G40}

\begin{abstract} We show that the two cuspidal unipotent characters of
a finite Chevalley group $E_7(q)$ have Schur index~$2$, provided that $q$ 
is an even power of a (sufficiently large) prime number $p$ such that 
$p\equiv 1 \bmod 4$. The proof uses a refinement of Kawanaka's generalized
Gelfand--Graev representations and some explicit computations with
the {\sf CHEVIE} computer algebra system.
\end{abstract}

\maketitle

\pagestyle{myheadings}
\markboth{Geck}{Schur indices of cuspidal unipotent characters of $E_7(q)$}

\section{Introduction} \label{MGsec0}

Throughout this paper, let $G$ be a simple algebraic group of adjoint 
type $E_7$. Assume that $G$ is defined over the finite field ${\F}_q$, 
with corresponding Frobenius map $F\colon G \rightarrow G$. There are 
precisely two cuspidal unipotent characters of $G^F$ denoted by $E_7[\pm 
\xi]$ where $\xi=\sqrt{-q}$; see the table in \cite{Ca2}, \S 13.9.  

The purpose of this paper is to determine the Schur index of $E_7[\pm \xi]$,
at least if the characteristic of ${\F}_q$ is large enough. Modulo this
condition on the characteristic, this completes the determination of the 
Schur indices of the unipotent characters of finite groups of Lie type; 
see \cite{LuRat}, \cite{mythom} and the references there.

By \cite{myert03}, Table~1, the character values of $E_7[\pm \xi]$ generate 
the field ${\Q}(\xi)$. Furthermore, by \cite{myert03}, Example~6.4, we 
already know that the Schur index is $1$ if $p\not\equiv 1 \bmod 4$ or if 
$q$ is not a square, where $p$ is the characteristic of $\F_q$. Thus, the
remaining task is to determine the Schur index when $q$ is a square and 
$p\equiv 1 \bmod 4$. 

\begin{thm} \label{thm1} Assume that $q$ is an even power of a (sufficiently
large) prime $p$ such that $p\equiv 1 \bmod 4$. Then the characters 
$E_7[\pm \xi]$ have Schur index~$2$.
\end{thm}

Here, $p$ is ``sufficiently large''  if Lusztig's results \cite{LuUS} on 
generalized Gelfand--Graev characters hold; it is conjectured that this
is the case if $p$ is good for $G$.

The idea of the proof is as follows. We have already seen in \cite{mythom},
\S 4, that $E_7[\pm \xi]$ occur with multiplicity~$1$ in a generalized 
Gelfand--Graev character $\Gamma_u$, where $u$ is a certain unipotent 
element in $G$. Here, we shall use a refinement of the construction of 
$\Gamma_u$ to show that, under the given assumptions on $p$ and $q$, the
characters $E_7[\pm \xi]$ occur with odd multiplicity in an induced character
which cannot be realized over $\Q(\xi)$. By standard arguments 
on Schur indices, this implies that $E_7[\pm \xi]$ cannot be realized  
over ${\Q}(\xi)$. At some stage, the proof relies on the fact that, in 
Lusztig's parametrization of the irreducible characters of $G^F$, the 
function $\Delta$ occurring in \cite{Lu4}, Main Theorem~4.23, takes value 
$-1$ on the labels corresponding to $E_7[\pm \xi]$.

Furthermore, we rely on some explicit computations in $G^F$. However, we 
shall only use computations with the root system and the irreducible 
characters of the Weyl group of $G$, for which the {\sf CHEVIE} system 
\cite{chv} is a convenient tool.

\section{Generalized Gelfand--Graev characters for type $E_7$} \label{MGsec1}

A short summary of the construction of generalized Gelfand--Graev characters
is given in \cite{mythom}, \S 2. Assume that $q$ is a power of a ``good''
prime $p\neq 2,3$.  Let $\Phi$ be the root system of $G$ with respect to a 
fixed maximally split torus $T$. Let $C$ be the unipotent class of $G$ whose
weighted 
Dynkin diagram $d\colon \Phi \rightarrow \Z$ is given in Table~\ref{tabe7}. 
(The notation in that table also defines a labelling of the simple roots in 
the root system of $G$.)  The class $C$ is the ``unipotent support'' of the 
two cuspidal unipotent characters of $G^F$; see \cite{mythom}, \S 4, and 
the references there.

\begin{table}[htbp] \caption{The weighted Dynkin diagram for the unipotent
support of the cuspidal unipotent characters in type $E_7$}
\label{tabe7}
\begin{center}
\begin{picture}(280,60)
\put( 20, 45){$E_7$}
\put( 60, 45){\circle*{8}}
\put( 60, 45){\line(1,0){40}}
\put(100, 45){\circle*{8}}
\put(100, 45){\line(1,0){40}}
\put(140, 45){\circle*{8}} 
\put(140, 45){\line(1,0){40}}
\put(180, 45){\circle*{8}}
\put(180, 45){\line(1,0){40}}
\put(220, 45){\circle*{8}}
\put(220, 45){\line(1,0){40}}
\put(260, 45){\circle*{8}}
\put(140, 45){\line(0,-1){40}}
\put(140,  5){\circle*{8}}
\put( 56, 31){$\alpha_1$}
\put( 57, 53){$1$}
\put( 96, 31){$\alpha_3$}
\put( 97, 53){$0$}
\put(126, 33){$\alpha_4$}
\put(137, 53){$1$}
\put(176, 31){$\alpha_5$}
\put(177, 53){$0$}
\put(216, 31){$\alpha_6$}
\put(217, 53){$1$}
\put(256, 31){$\alpha_7$}
\put(257, 53){$0$}
\put(149,  2){$\alpha_2$}
\put(125,  1){$0$}
\end{picture}
\end{center}
\end{table}

Given the weight function $d\colon \Phi \rightarrow \Z$ specified by the 
diagram in Table~\ref{tabe7}, we define unipotent subgroups
\[ U_{d,2}:=\prod_{\atop{\alpha\in \Phi^+}{d(\alpha)\geq 2}} X_\alpha
\qquad \mbox{and}\qquad U_{d,1}:=\prod_{\atop{\alpha\in \Phi^+}{d(\alpha)
\geq 1}} X_\alpha,\]
where $X_\alpha$ is the root subgroup in $G$ corresponding to the root
$\alpha$. (It is understood that the products are taken in some fixed order.)
The generalized Gelfand--Graev character associated with an element in
$C^F$ is obtained by inducing a certain linear character from $U_{d,2}^F$.
We have $C_G(u)/C_G(u)^\circ \cong {\Z}/2{\Z}$ for $u\in C$. Thus, $C^F$ 
splits into two classes in the finite group $G^F$. By Mizuno \cite{Miz}, 
Lemma~28, representatives of these two $G^F$-classes are given by 
\begin{align*}
y_{74}&= x_{20}(1)x_{21}(1)x_{23}(1)x_{28}(1)x_{31}(1),\\
y_{75}&= x_{20}(1)x_{21}(1)x_{28}(1)x_{24}(1)x_{23}(1)x_{25}(1)x_{36}
(\zeta),
\end{align*}
where $\zeta$ is a generator for the multiplicative group of $\F_q$ and 
where the subscripts correspond to the following roots in $\Phi^+$:
\begin{alignat*}{2}
20&: \alpha_1+\alpha_2+\alpha_3+\alpha_4,&\qquad 
21&: \alpha_1+\alpha_3+\alpha_4+\alpha_5,\\
23&: \alpha_2+\alpha_4+\alpha_5+\alpha_6,&\qquad
24&: \alpha_3+\alpha_4+ \alpha_5+\alpha_6,\\
25&: \alpha_4+ \alpha_5+\alpha_6+\alpha_7,&\qquad
28&: \alpha_2+\alpha_3+2\alpha_4+\alpha_5,\\
31&: \alpha_3+\alpha_4+ \alpha_5+\alpha_6+\alpha_7,&\qquad
36&: \alpha_2+\alpha_3+\alpha_4+ \alpha_5+\alpha_6+\alpha_7.
\end{alignat*}
(Attention: Here, we use the labelling of the roots as given by  the
{\sf CHEVIE} system \cite{chv}, which is slightly different from that of
Mizuno.) We note that both $y_{74}$ and $y_{75}$ lie in $C \cap U_{d,2}^F$. 
Now let us fix 
\[ u\in \{y_{74},y_{75}\}\subseteq C\cap U_{d,2}^F;\]
the above expressions show that 
\[ u=\prod_{\atop{\alpha \in \Phi^+}{d(\alpha)=2}} x_\alpha(
\eta_\alpha)\qquad\mbox{where $\eta_\alpha\in \F_q$}.\]
Then we define a linear character $\varphi_u\colon U_{d,2}^F \rightarrow 
\C^\times$ by the formula 
\[ \varphi_u\Bigl(\prod_{\atop{\alpha\in \Phi^+}{d(\alpha)\geq 2}}x_\alpha
(\xi_\alpha)\Bigr)=\chi\Bigl(\sum_{\atop{\alpha\in \Phi^+}{d(\alpha)=2}}
c_\alpha\,\eta_{\alpha}\, \xi_{\alpha}\Bigr) \qquad \mbox{for all 
$\xi_\alpha\in \F_q$},\]
where $c_\alpha\in \F_q$ are certain fixed constants (independent of
the $\eta_\alpha$ and $\xi_\alpha$) and where $\chi\colon {\F}_q^+
\rightarrow\C^\times$ is a fixed non-trivial character of the additive 
group of ${\F}_q$; see \cite{mythom}, Definition~2.1, for more details. 
It will actually be convenient to choose $\chi$ in the following special
way. Let $\chi_0\colon {\F}_p^+\rightarrow\C^\times$ be a fixed 
non-trivial 
character of the additive group of ${\F}_p$. Then we take $\chi$ to be
\[ \chi:=\chi_0\circ \mbox{Tr}_{{\F}_q/{\F}_p}\]
where $\mbox{Tr}_{{\F}_q/{\F}_p}\colon {\F}_q^+\rightarrow {\F}_p^+$ is the
trace map. Now we have 
\[ \Ind_{U_{d,2}^F}^{G^F}\bigl(\varphi_u\bigr)=[U_{d,1}^F:U_{d,2}^F]^{1/2}
\cdot \Gamma_u,\]
where $\Gamma_u$ is the generalized Gelfand--Graev character associated
with~$u$. We have seen in \cite{mythom}, Corollary~4.3, that 
\[ \big\langle E_7[\pm \xi],\Gamma_u\big\rangle_{G^F}=1 \qquad 
\mbox{for suitable $u\in \{y_{74}, y_{75}\}$}.\]
(Here, and throughout the paper, we denote by $\langle\;,\;\rangle_A$
the standard inner product on the character ring of a finite group $A$.)

We will now refine the construction of $\Gamma_u$. The strategy for doing 
this has already been outlined in \cite{mythom}, \S 4. For this purpose,
we shall assume from now on that
\begin{center}
{\em $q$ is an even power of $p$}.
\end{center}
Since $G$ is simple of adjoint type, we have an ${\F}_q$-isomorphism 
\[ h\colon \underbrace{k^\times \times \cdots\times  k^\times}_{\text{$7$ 
factors}} \rightarrow T, \qquad (x_1,\ldots,x_7)\mapsto h(x_1,\ldots,x_7),\]
such that $\alpha_i(h(x_1,\ldots,x_7))=x_i$ for $1\leq i \leq 7$. In 
particular, we have $T^F=\{h(x_1,\ldots,x_7) \mid x_i\in {\F}_q^\times\}$.
We shall set 
\[ t:=h(\nu^{1/2},1,1,\nu^{1/2},1, \nu^{1/2},1)\in T^F\] 
as in the proof of \cite{mythom}, Lemma~4.1, where $\nu$ is a generator 
for the multiplicative group of $\F_p\subset \F_q$ and $\nu^{1/2}$ is a 
square root of $\nu$ in $\F_q$. (The square root exists since $q$ is an
even power of $p$.) Then $t$ has the property that $\alpha(t)=\nu$ for 
all roots $\alpha$ involved in the expressions for $y_{74}$ or $y_{75}$
as products of root subgroup elements; furthermore, we have $\alpha(t)=1$ 
for all roots $\alpha$ such that $d(\alpha)=0$. The element $t$ has order 
$2(p-1)$ and $H:=\langle t\rangle$ normalizes $U_{d,2}$. We set 
\[ s_1:=h(-1,1,1,-1,1,-1,1)=t^{p-1} \in T^F.\]
Note that $\alpha(s_1)=1$ for all roots $\alpha\in \Phi^+$ which are 
involved in the expressions of $y_{74}$ and $y_{75}$ as products of root 
subgroup elements. Thus, $s_1$ fixes the character $\varphi_u$ and so we 
can extend $\varphi_u$ to $U_{d,2}^F.\langle s_1\rangle$. Actually, there 
are two such extensions which we denote by $\tilde{\varphi}_u$ and 
$\tilde{\varphi}_u'$. Their values are determined by 
\[ \tilde{\varphi}_{u}(xs_1)=\varphi_u(x) \qquad \mbox{and}
\qquad \tilde{\varphi}_{u}'(xs_1)=-\varphi_u(x) \qquad \mbox{for all 
$x\in U_{d,2}^F$}.\]

\begin{defn} \label{def1} Let $u\in \{y_{74},y_{75}\}$. Then we set 
\[ \psi_u:=\Ind_{U_{d,2}^F.\langle s_1\rangle}^{U_{d,2}^F.H}\bigl(
\tilde{\varphi}_u\bigr) \qquad \mbox{and}\qquad 
\psi_u':=\Ind_{U_{d,2}^F.\langle s_1\rangle}^{U_{d,2}^F.H}\bigl(
\tilde{\varphi}_u'\bigr).\]
Thus, we have 
\[\Ind_{U_{d,2}^F}^{U_{d,2}^F.H}\bigl(\varphi_u\bigr)=\psi_u+\psi_u'
\qquad \mbox{and}\qquad [U_{d,1}^F:U_{d,2}^F]^{1/2}\cdot \Gamma_u=
\tilde{\Gamma}_u+ \tilde{\Gamma}_u',\]
where 
\begin{align*}
\tilde{\Gamma}_u&:=\Ind_{U_{d,2}^F.\langle s_1\rangle}^{G^F}
\bigl(\tilde{\varphi}_u\bigr)=\Ind_{U_{d,2}^F.H}^{G^F}\bigl(\psi_u\bigr),\\
\tilde{\Gamma}_u'&:=\Ind_{U_{d,2}^F.\langle s_1\rangle}^{G^F}
\bigl(\tilde{\varphi}_u'\bigr)=\Ind_{U_{d,2}^F.H}^{G^F}\bigl(\psi_u'\bigr).
\end{align*}
The following result provides some crucial information concerning 
$\psi_u$ and $\psi_u'$.
\end{defn}

\begin{prop} \label{lem1} Recall that $q$ is an even power of $p$. Then, 
with the above notation, the following hold.
\begin{itemize}
\item[(a)] Both $\psi_u$ and $\psi_u'$ are irreducible characters of
$U_{d,2}^F.H$.
\item[(b)] $\psi_u$ can be realized over $\Q$.
\item[(c)] $\psi_u'$ is rational-valued but cannot be realized over $\Q$.
In fact, $\psi_u'$ has non-trivial local Schur indices at $\infty$ and at
the prime $p$.
\end{itemize}
\end{prop}

\begin{proof} (See also the argument of Ohmori \cite{Ohm2}, p.~154.) Let
\begin{equation*} 
x:=\prod_{\atop{\alpha\in \Phi^+}{d(\alpha)\geq 2}}x_\alpha
(\xi_\alpha)\in U_{d,2}^F\qquad\mbox{and}\qquad \gamma_x:=
\sum_{\atop{\alpha\in \Phi^+}{d(\alpha)=2}} c_\alpha\eta_\alpha\xi_\alpha,
\tag{1}
\end{equation*}
where $\xi_\alpha\in {\F}_q$. Then, as in the proof of \cite{mythom}, 
Proposition~2.3, we have 
\[\varphi_u(t^ixt^{-i})=\chi(\nu^{i}\gamma_x)\qquad\mbox{for $1\leq i
\leq 2(p-1)$}.\]
In particular, this implies $\mbox{Stab}_H(\varphi_u)=\langle s_1
\rangle$. Hence, by Clifford theory, the induced character 
\[\Ind_{U_{d, 2}^F}^{U_{d,2}^F.H}\bigl(\varphi_u\bigr)=\psi_u+\psi_u'\]
has inner product~$2$. Thus, we $\psi_u$ and $\psi_u'$ must be irreducible, 
proving~(a).

Next we prove (b). Using Mackey's formula and relation (1), we 
have that 
\begin{align*}
\Ind_{U_{d,2}^F}^{U_{d,2}^F.H}&\bigl(\varphi_u\bigr)(x)=\sum_{i=1}^{2(p-1)} 
\varphi_u(t^ixt^{-i})=\sum_{i=1}^{2(p-1)}\chi(\nu^{i}\gamma_x)\\ &=
\sum_{i=1}^{2(p-1)} \chi_0\bigl(\nu^i\, \mbox{Tr}_{{\F}_q/{\F}_p}(\gamma_x)
\bigr)= \left\{\begin{array}{cl} 2(p-1) & \mbox{if $\mbox{Tr}_{{\F}_q/
{\F}_p}(\gamma_x)=0$},\\ -2 & \mbox{if $\mbox{Tr}_{{\F}_q/{\F}_p}(\gamma_x)
\neq 0$}.\end{array}\right.
\end{align*}
In particular, this shows that the values are rational integers. 
Thus, $\psi_u+\psi_u'$ is rational-valued. Now assume, if possible, 
that $\psi_u$ is not rational-valued. Then the characters $\psi_u$ 
and $\psi_u'$ must be algebraically conjugate. Consequently, 
$\psi_u$ and $\psi_u'$ occur with the same multiplicity in every 
rational-valued character. Now, by the Mackey formula and Frobenius 
reciprocity, we have 
\begin{align*}
\Big\langle \psi_u',\Ind_H^{U_{d,2}^F.H}\bigl({\bf 1}_{H}\bigr)\Big 
\rangle_{U_{d,2}^F.H}&=\Big\langle \Ind_{U_{d,2}^F.\langle s_1
\rangle}^{U_{d,2}^F.H} \bigl(\tilde{\varphi}_u' \bigr),\Ind_H^{U_{d,2}^F.H}
\bigl({\bf 1}_{H}\bigr) \Big \rangle_{U_{d,2}^F.H}\\ &=\big\langle
\Res_{\langle s_1 \rangle}^{U_{d,2}^F}\bigl( \tilde{\varphi}_u' \bigr),
{\bf 1}_{\langle s_1\rangle}\big \rangle_{U_{d,2}^F}\\ &=0,
\end{align*}
since $\tilde{\varphi}_u'(s_1)=-1$. (Here, the symbol ${\bf 1}$ stands
for the unit character.) By a similar argument, since $\tilde{\varphi}_u
(s_1)=1$,  we also have
\[\Big\langle \psi_u,\Ind_H^{U_{d,2}^F.H}\bigl({\bf 1}_{H}\bigr)\Big 
\rangle_{U_{d,2}^F.H}=1.\]
Thus, $\psi_u$ and $\psi_u'$ do not occur with the same multiplicity in
some rational-valued character, a contradiction. Thus, our assumption was 
wrong and so both $\psi_u$ and $\psi_u'$ are rational-valued. But then the 
above multiplicity~$1$ formula implies that $\psi_u$ can be realized over
$\Q$, by a standard argument concerning Schur indices (see Isaacs 
\cite{Isa}, Corollary~10.2). 

Finally, we prove (c). We begin by showing that the local Schur index
at $\infty$ is non-trivial. In other words, we must show that $\psi_u'$
cannot be realized over $\R$. For this purpose, by a well-known criterion 
due to Frobenius and Schur (see Isaacs \cite{Isa}, Chapter~4), it is enough 
to show that
\[\frac{1}{|U_{d,2}^F.H|}\sum_{g\in U_{d,2}^F.H} \psi_u'(g^2)=-1.\]
Now, in order to evaluate the above sum, we note that
\[\frac{1}{|U_{d,2}^F.H|} \sum_{g\in U_{d,2}^F.H} \psi_u(g^2)=1,\]
since $\psi_u$ can be realized over $\Q$. Thus, it will be enough to show
that 
\[\frac{1}{|U_{d,2}^F.H|} \sum_{g\in U_{d,2}^F.H} 
\Ind_{U_{d,2}}^{U_{d,2}^F.H}\bigl(\varphi_u \bigr) (g^2)=
\frac{1}{|U_{d,2}^F.H|}\sum_{g\in U_{d,2}^F.H} (\psi_u+\psi_u')(g^2)=0.\]
Let $g\in U_{d,2}^F.H$ and write $g=xh$ where $x\in U_{d,2}^F$ and $h\in H$.
Now the value of the above induced character on $g^2$ is zero unless
$g^2\in U_{d,2}^F$. Thus, we only need to consider elements $g=xh$ where
$h=1$ or $h=s_1$. So we must show that 
\[\sum_{x\in U_{d,2}^F} \Ind_{U_{d,2}}^{U_{d,2}^F.H}\bigl(\varphi_u\bigr)
(x^2)+\sum_{x\in U_{d,2}^F}\Ind_{U_{d,2}}^{U_{d,2}^F.H}\bigl(\varphi_u\bigr) 
(xs_1xs_1)=0.\]
Now, since $U_{d,2}^F$ has odd order, the map $x\mapsto x^2$ defines a
bijection of $U_{d,2}^F$ onto itself. Hence the first sum evaluates to
\[\sum_{x\in U_{d,2}^F} \Ind_{U_{d,2}}^{U_{d,2}^F.H}\bigl(\varphi_u\bigr)
(x^2)=|U_{d,2}^F.H|\cdot \Big\langle \Ind_{U_{d,2}}^{U_{d,2}^F.H}
\bigl(\varphi_u\bigr), {\bf 1}_{U_{d,2}^F.H}\Big\rangle_{U_{d,2}^F.H}=0.\]
Now consider the second sum. For this purpose, we note that $\alpha(s_1)=1$ 
for all roots $\alpha\in \Phi^+$ which are involved in the expressions of 
$y_{74}$ and $y_{75}$ as products of root subgroup elements. Thus, if $x$ 
and $\gamma_x$ are as in (1), then we have 
\[\gamma_{(xs_1)^2}=\sum_{\atop{\alpha\in \Phi^+}{d(\alpha)=2}} 
c_\alpha\eta_\alpha(\alpha(s_1)+1)\xi_\alpha =
2\gamma_x=\gamma_{x^2}.\]
Using once more Mackey's formula as at the beginning of this proof,
we see  that 
\begin{align*}
\Ind_{U_{d,2}}^{U_{d,2}^F.H}\bigl(\varphi_u\bigr)(x^2) &=
\sum_{i=1}^{2(p-1)}\chi(\nu^{i}\gamma_{x^2})= \sum_{i=1}^{2(p-1)}
\chi(\nu^{i}\gamma_{(xs_1)^2})\\&= \Ind_{U_{d,2}}^{U_{d,2}^F.H}
\bigl(\varphi_u\bigr)(xs_1xs_1)
\end{align*}
for all $x\in U_{d,2}^F$. Consequently, the second sum also equals~$0$.
Thus, we have shown that $\psi_u'$ cannot be realized over $\R$. We shall 
now use some general properties of Schur indices; see Feit \cite{Feit1},
\S 2, for references. First, since $\psi_u'$ is rational-valued
but $\psi_u'$ cannot be realized over $\R$, the Schur index of $\psi_u'$ 
is~$2$ (by the Brauer--Speiser theorem; see \cite{Feit1}, 2.4). Furthermore,
there exists at least one prime number $\ell$ such that the $\ell$-local 
Schur index of $\psi_u'$ is~$2$ (by the Hasse sum formula; see \cite{Feit1}, 
2.15). Thus, it will be enough to show that the $\ell$-local Schur index of 
$\psi_u'$ is $1$, for every prime $\ell\neq p$. Let $\ell$ be such a prime. 
If $\ell\neq 2$, then $\psi_u'$ is a character of 
$\ell$-defect $0$ of $U_{d,2}^F.H$. So the $\ell$-local Schur index is 
$1$ by \cite{Feit1}, 2.10. Finally, if $\ell=2$, then $\psi_u'$ is a 
character of $2$-defect~$1$ and, hence, lies in a block with a cyclic 
defect group of order~$2$. Consequently, that block contains only two 
irreducible characters and so $\psi_u'$ remains irreducible as a 
$2$-modular Brauer character. This implies again that the local Schur
index is $1$; see \cite{Feit1}, 2.10.
\end{proof}

\section{A subgroup of type $D_6 \times A_1$} \label{MGsec2}
Our next aim is to compute the multiplicity of $E_7[\pm \xi]$ in 
$\tilde{\Gamma}_u$ and $\tilde{\Gamma}_u'$; see Definition~\ref{def1}. 
We already know that the multiplicity of $E_7[\pm \xi]$ in the sum 
$\tilde{\Gamma}_u+ \tilde{\Gamma}_u'$ equals $[U_{d,1}^F:U_{d,2}^F]^{1/2}$, 
for suitable $u\in \{y_{74},y_{75}\}$. We shall now try to compute the 
multiplicity in the difference $\tilde{\Gamma}_u-\tilde{\Gamma}_u'$. For
this purpose, we take a closer look at the semisimple element $s_1$ and its 
centralizer. Let 
\[ G_1:=\langle T,X_\alpha \mid \alpha \in \Phi_1\rangle\qquad\mbox{where}
\qquad \Phi_1:=\{\alpha\in \Phi\mid \alpha(s_1)=1\}.\]
Using the {\sf CHEVIE} function {\tt ReflectionSubgroup}, we check that 
the root system $\Phi_1$ has type $D_6 \times A_1$; a system of simple 
roots in $\Phi_1$ is given by 
\[ \Pi_1=\{\alpha_2,\alpha_3,\alpha_5,\alpha_7,\alpha_{14},\alpha_{18},
\alpha_{28}\}\]
where 
\[ \alpha_{14}:=\alpha_1{+}\alpha_3{+}\alpha_4,\quad    
\alpha_{18}:=\alpha_4{+}\alpha_5{+}\alpha_6,\quad       
\alpha_{28}:= \alpha_2{+}\alpha_3{+}2\alpha_4{+}\alpha_5.\]
(Here, the numbering of the roots is the same as that given by {\sf CHEVIE}.)
The corresponding Dynkin diagram and the restriction of the weight 
function $d$ to $\Pi_1$ are given in Table~\ref{tabd6a1}. Furthermore, one 
can check, using {\sf CHEVIE} (for example), that
\[ N_W(W_1)=\{w\in W \mid w(\Phi_1)\subseteq \Phi_1\}=W_1\]
where $W_1:=\langle w_\alpha \mid \alpha \in \Phi_1\rangle \subset W$ is
the Weyl group of $G_1$ (and where we denote by $w_\beta$ the reflection 
with root $\beta$, for any root $\beta\in \Phi$).  

\begin{lem} \label{lem3} We have $C_G(s_1)=G_1$; in particular, $C_G(s_1)$ 
is connected. 
\end{lem}

\begin{proof} By Carter \cite{Ca2}, \S~3.5, we have $C_G(s_1)^\circ=
G_1$. Hence, $G_1$ is a normal subgroup in $C_G(s_1)$. So it is enough to 
show that $N_G(G_1)=G_1$. Let $g\in N_G(G_1)$. Then $gTg^{-1}$ is a maximal 
torus in $G_1$ and so there exists some $g_1\in G_1$ such that $gTg^{-1}=
g_1Tg_1^{-1}$. Thus, we have $g_1^{-1}g \in N_G(T)$ and so $g\in G_1.N_G(T)$.
Hence, we may assume without loss of generality that $g\in N_G(T)\cap N_G
(G_1)$. Now, for any $g\in N_G(T)\cap N_G(G_1)$ and any $\alpha\in \Phi_1$, 
we have $gX_\alpha g^{-1}=X_{w(\alpha)} \subseteq G_1$, where $w$ is the image 
of $g$ in $W=N_G(T)/T$. Thus, we have $w(\Phi_1)\subseteq \Phi_1$ and so
$w\in W_1$ (see the above remarks). This implies $g\in G_1$, as required. 
\end{proof}

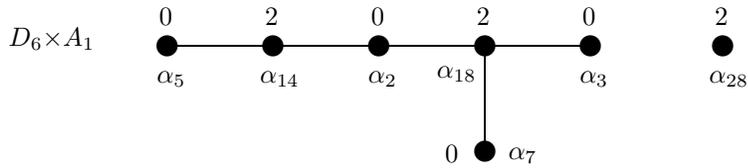
\begin{table}[htbp] \caption{The restriction of $d$ to the subsystem of
type $D_6 \times A_1$}
\label{tabd6a1}
\begin{center}
\begin{picture}(330,60)
\put(  0, 45){$D_6{\times} A_1$}
\put( 60, 45){\circle*{8}}
\put( 60, 45){\line(1,0){40}}
\put(100, 45){\circle*{8}}
\put(100, 45){\line(1,0){40}}
\put(140, 45){\circle*{8}} 
\put(140, 45){\line(1,0){40}}
\put(180, 45){\circle*{8}}
\put(180, 45){\line(1,0){40}}
\put(220, 45){\circle*{8}}
\put(270, 45){\circle*{8}}
\put(180, 45){\line(0,-1){40}}
\put(180,  5){\circle*{8}}
\put( 56, 31){$\alpha_5$}
\put( 57, 53){$0$}
\put( 95, 31){$\alpha_{14}$}
\put( 97, 53){$2$}
\put(136, 31){$\alpha_2$}
\put(137, 53){$0$}
\put(162, 33){$\alpha_{18}$}
\put(177, 53){$2$}
\put(216, 31){$\alpha_3$}
\put(217, 53){$0$}
\put(265, 31){$\alpha_{28}$}
\put(267, 53){$2$}
\put(189,  2){$\alpha_7$}
\put(165,  1){$0$}
\end{picture}
\end{center}
\end{table}


Let $C_1$ be the conjugacy class of $y_{74}$ in $G_1$ and denote by $d_1
\colon \Phi_1\rightarrow \Z$ the corresponding weighted Dynkin diagram. 
Using the identification results in \cite{Ca1}, Theorem~11.3.2, it is 
straightforward to check that, under the natural matrix representation of 
a group of type $D_6 \times A_1$, the elements $y_{74}$ and $y_{75}$ 
correspond to matrices with Jordan blocks of size $1,1,2,5,5$ (where 
the block of size $2$ comes from the $A_1$-factor). Hence, using 
\cite{Ca2}, \S 13.1, we see that $d_1$ is given by the restriction of $d$
to $\Phi_1$, as specified in Table~\ref{tabd6a1}. Furthermore, we notice
that the above roots can all be written as sums of roots in $\Pi_1$.  
Thus, we have 
\[ y_{74}, y_{75} \in C_1\cap U_{d_1,2}^F,\]
where $U_{d_1,2}$ is the unipotent subgroup of $G_1$ defined with respect
to $d_1$.

\begin{lem} \label{lem2} Let $u\in \{y_{74},y_{75}\}$. Then we have
$\dim \fB_u^1=4$ (where $\fB_u^1$ denotes the variety of Borel subgroups
of $G_1$ containing $u$) and 
\[ C_{G_1}(u)/C_{G_1}(u)^\circ \cong C_G(u)/C_G(u)^\circ \cong {\Z}/2{\Z}.\] 
\end{lem}

\begin{proof} Let $u:=y_{74}$. The formula for $\dim \fB_u^1$ follows from 
\cite{Ca2}, \S 13.1.  To prove the remaining statements, we note that 
\[ s_1\in S:=\{h(x,x^{-2},x^{-2},x^3,x^{-2},x,1) \mid x\in k^\times \}
\subseteq C_{G_1}(u).\]
Furthermore, one checks that $Z(G_1)=\{t\in T\mid \alpha(t)=1\mbox{ for all
$\alpha\in \Phi_1$}\}=\langle s_1\rangle$. Thus, since $S$ is connected, we 
have $Z(G_1)\subseteq C_{G_1}(u)^\circ$.

Now let $\pi \colon G_1\rightarrow H_1$ be the adjoint quotient of $G_1$, 
where $H_1$ is a semisimple group of adjoint type $D_6\times A_1$. Let 
$\bar{u}$ be the image of $u$ in $H_1$. Then, by Carter \cite{Ca2}, \S 13.1,
we know that $C_{H_1}(\bar{u})/C_{H_1}(\bar{u})^\circ \cong 
{\Z}/2{\Z}$. Furthermore, $\pi$ induces a surjective homomorphism 
\[ C_{G_1}(u)/C_{G_1}(u)^\circ\twoheadrightarrow C_{H_1}(\bar{u})/C_{H_1}
(\bar{u})^\circ \cong {\Z}/2{\Z}\]
with kernel given by the image of $Z(G_1)$ in $C_{G_1}(u)/C_{G_1}(u)^\circ$.
Since $Z(G_1)\subseteq C_{G_1}(u)^\circ$, that image is trivial and so the
above surjective map is also injective. 
\end{proof}

\begin{prop} \label{prop1} Let $u \in \{y_{74},y_{75}\}\subseteq C\cap 
U_{d,2}^F$.  Then, as we already noted, we have $u\in C_1\cap U_{d_1,2}^F$ 
and so the corresponding generalized Gelfand--Graev character $\Gamma_{u}^1$ 
of $G_1^F$ is well-defined. We have 
\[ \tilde{\Gamma}_u(ys_1)-\tilde{\Gamma}_u'(ys_1)=\Gamma_{u}^1(y)
\qquad \mbox{for all $y\in G_1^F$ unipotent}.\]
\end{prop} 

\begin{proof} By the Mackey formula, we have
\begin{align*}
\tilde{\Gamma}_u(ys_1) &= \Res_{G_1^F}^{G^F}\bigl(\tilde{\Gamma}_u
\bigr)(ys_1)=\Res_{G_1^F}^{G^F}\Bigl(\Ind_{U_{d,2}^F.\langle s_1 
\rangle}^{G^F} \bigl(\tilde{\varphi}_u\bigr)\Bigr)(ys_1)\\&=\sum_z
\Ind_{(U_{d,2}^F.\langle s_1\rangle)^z \cap G_1^F}^{G_1^F}\Bigl(
\Res_{(U_{d,2}^F.\langle s_1\rangle)^z \cap G_1^F}^{(U_{d,2}^F.\langle 
s_1\rangle)^z}\bigl(\tilde{\varphi}_u^z\bigr)\Bigr)(ys_1),
\end{align*}
where $z$ runs over a set of representatives of the $(U_{d,2}^F.\langle s_1
\rangle,G_1^F)$-double cosets of $G^F$. Let us fix such a double coset 
representative, $z$ say.  Assume that the value at $ys_1$ of the 
corresponding induced character in the above sum is non-zero. Then $ys_1$ 
must be $G_1^F$-conjugate to an element in the subgroup $(U_{d,2}^F.\langle 
s_1\rangle)^z \cap G_1^F$. Consequently, $s_1$ must be $G_1^F$-conjugate to 
an element in that subgroup. Since $\langle s_1\rangle$ is a Sylow
$2$-subgroup of $U_{d,2}^F.\langle s_1\rangle$, we conclude that all
elements of order~$2$ in $U_{d,2}^F.\langle s_1\rangle$ are of the form
$xs_1x^{-1}$ where $x\in U_{d,2}^F$. Thus, we have $c^{-1}s_1c=z^{-1}xs_1
x^{-1}z$ for some $c\in G_1^F$ and some $x\in U_{d,2}^F$. Consequently, 
$x^{-1}zc^{-1} \in C_G(s_1)^F=G_1^F$ and so $z\in xG_1^Fc\in 
U_{d,2}^F.G_1^F$. Thus, $z$ represents the trivial double coset and so
we can take $z=1$. Using the fact that 
\[ U_{d,2}^F.\langle s_1\rangle\cap G_1^F=U_{d_1,2}^F \times \langle s_1
\rangle\]
(where $U_{d_1,2}\subseteq G_1$ is the unipotent subgroup defined with
respect to the weighted Dynkin diagram $d_1\colon \Phi_1\rightarrow \Z$) 
we find that 
\[\tilde{\Gamma}_u(ys_1) = \Ind_{U_{d_1,2}^F\times \langle s_1
\rangle}^{G_1^F}\bigl(\varphi_u^1 \boxtimes {\bf 1}_{\langle s_1\rangle}
\bigr)(ys_1)\]
where $\varphi_u^1$ denotes the restriction of $\varphi_u$ to 
$U_{d_1,2}^F$. Since $s_1$ is in the center of $G_1$, it is readily 
checked that 
\[\tilde{\Gamma}_u(ys_1) = \frac{1}{2}\,\tilde{\varphi}_u(s_1)\,
\Ind_{U_{d_1,2}^F}^{G_1^F}\bigl(\varphi_u^1\bigr)(y)=\frac{1}{2}\,
\Ind_{U_{d_1,2}^F}^{G_1^F}\bigl(\varphi_u^1\bigr)(y).\]
By a completely analogous argument, we also obtain that 
\[\tilde{\Gamma}_u'(ys_1) = \frac{1}{2}\,\tilde{\varphi}_u'(s_1)\,
\Ind_{U_{d_1,2}^F}^{G_1^F}\bigl(\varphi_u^1\bigr)(y)=-\frac{1}{2}
\Ind_{U_{d_1,2}^F}^{G_1^F}\bigl(\varphi_u^1\bigr)(y).\]
Thus, it remains to check that 
\[ \Gamma_{u}^1=\Ind_{U_{d_1,2}^F}^{G_1^F}\bigl(\varphi_u^1\bigr).\]
For this purpose, we must show that $\varphi_u^1$ indeed is the linear 
character of $U_{d_1,2}^F$ required in the definition of $\Gamma_u^1$.
Now, the definition of $\Gamma_u^1$ requires the choice of a non-degenerate
bilinear form and of an oppisition automorphism on the Lie algebra of $G_1$.
However, the Lie algebra of $G_1$ is naturally contained in the Lie
algebra of $G$, with compatible Cartan decompositions. Thus, the chosen 
bilinear form and the chosen opposition automorphism restrict to the Lie 
algebra of $G_1$, and this implies that $\varphi_u^1$ is the required linear
character of $U_{d_1,2}^F$. 
\end{proof}

A formula of this kind has been stated (without proof) by Kawanaka in
\cite{Kaw3}, Lemma~2.3.5; see also the Ph. D. thesis of Wings \cite{Wings}, 
\S 3.2.1. 

\begin{rem} \label{rem1} Let $g\in G^F$ and write $g=g_sg_u=g_ug_s$ where
$g_s\in G^F$ is semisimple and $g_u\in G^F$ is unipotent. Assume that
$g_s$ is not conjugate to $s_1$ in $G^F$. Then we have 
\[ (\tilde{\Gamma}_u-\tilde{\Gamma}_u')(g)=0.\]
Indeed, if the value is non-zero, then $g$ must be $G^F$-conjugate to
an element in $U_{d,2}^F.\langle s_1\rangle$. But then $g_s$ will also
be $G^F$-conjugate to an element in that subgroup. Using a Sylow
argument as in the above proof, we see that either $g_s=1$ or $g_s$ is 
$G^F$-conjugate to $s_1$, as claimed. Furthermore, if $g_s=1$, then it
is readily checked that $\tilde{\Gamma}_u(g)=\tilde{\Gamma}_u'(g)$.
\end{rem}

Thus, in order to compute the scalar product of $E_7[\pm \xi]$ with
$\tilde{\Gamma}_u-\tilde{\Gamma}_u'$, it will be enough to know the
values of $E_7[\pm \xi]$ on elements of the form $ys_1$ where $y\in G_1^F$
is unipotent. Furthermore, since $E_7[\xi]$ and $E_7[-\xi]$ are complex
conjugate and since $\tilde{\Gamma}_u$ and $\tilde{\Gamma}_u'$ are
rational-valued, it will actually be enough to consider the sum 
$E_7[\xi]+E_7[-\xi]$. Now, by Lusztig \cite{Lu4}, Main Theorem~4.23, we have
\[ E_7[\xi]+E_7[-\xi]=R_{512_a}-R_{512_a'}.\]
(Note that the function $\Delta$ occurring in \cite{Lu4}, 4.23, takes
value $-1$ on the labels corresponding to the characters $E_7[\pm \xi]$.)
Here, $512_a$, $512_a'$ are the two irreducible characters of $W$ of
degree~$512$ and $R_{512_a}$, 
$R_{512_a'}$ are the corresponding ``almost characters'', as defined by 
Lusztig \cite{Lu4}, (3.7). For any $\phi \in \Irr(W)$, we have 
\[ R_\phi:=\frac{1}{|W|} \sum_{w\in W} \phi(w)\, R_{T_w,1};\]
here, $T_w\subseteq G$ is an $F$-stable maximal torus  obtained from
$T$ by twisting with~$w$  and $R_{T_w,1}$ is the Deligne--Lusztig
generalized character associated with the trivial character of $T_w^F$.
Similarly, for any $\psi\in \Irr(W_1)$, we denote by $R_\psi^1$ the 
corresponding almost character of $G_1^F$.

\begin{lem} \label{lem4} Let $\phi \in \Irr(W)$ and write
\[ \Res_{W_1}^W\bigl(\phi\bigr)=\sum_{\psi \in \Irr(W_1)} m(\phi,\psi)\,
\psi \qquad \mbox{where $m(\phi,\psi)\in {\Z}_{\geq 0}$}.\]
Let $y\in G_1^F$ be a unipotent element.  Then we have 
\[ R_\phi(ys_1)=\sum_{\psi\in \Irr(W_1)} m(\phi,\psi)\, R_{\psi}^1(y).\]
\end{lem}

\begin{proof} The character formula for $R_{T_w,1}$ (see \cite{Ca2}, 
Theorem~7.2.8) shows that
\[ R_{T_w,1}(ys_1)=\frac{|C_W(w)|}{|W_1|}\sum_{\atop{w_1\in W_1}{w\sim w_1}}
\, R_{T_{w_1},1}^1(y)\]
where the relation $\sim$ means conjugacy in $W$. (Here, $R_{T_{w_1},1}^1$ 
denotes a Deligne--Lusztig generalized character of $G_1^F$.) Thus, we have 
\begin{align*}
R_\phi(ys_1) &= \frac{1}{|W|}\sum_{\atop{w\in W,w_1\in W_1}{w\sim w_1}}
\frac{|C_W(w)|}{|W_1|}\, \phi(w)\, R_{T_{w_1},1}^1(y)\\
&= \frac{1}{|W_1|}\sum_{w_1\in W_1} \Bigl(\frac{1}{|W|}
\sum_{\atop{w\in W}{w\sim w_1}} |C_W(w)|\,\phi(w)\Bigr) R_{T_{w_1},1}^1(y)
\end{align*}
Now, we have $\phi(w)=\phi(w_1)$ and $|C_W(w)|=|C_W(w_1)|$ for all $w_1
\in W_1$ such that $w\sim w_1$. Thus, we have
\[\frac{1}{|W|} \sum_{\atop{w\in W}{w\sim w_1}} |C_W(w)|\,\phi(w)
=\frac{|C_W(w_1)|}{|W|} \phi(w_1)\sum_{\atop{w\in W}{w\sim w_1}} 1=
\phi(w_1).\]
Writing $\phi(w_1)=\sum_\psi m(\phi,\psi)\psi(w_1)$, we obtain the
desired expression.
\end{proof}

\begin{cor} \label{cor1} With the notation of Proposition~\ref{prop1}
and Lemma~\ref{lem4}, we have
\[ \big\langle R_\phi,\tilde{\Gamma}_u-\tilde{\Gamma}_u'\big\rangle_{G^F}=
\sum_{\psi\in \Irr(W_1)} m(\phi,\psi)\,\big\langle R_\psi,\Gamma_u^1
\big\rangle_{G_1^F},\]
for any $\phi \in \Irr(W)$ and $u\in \{y_{74},y_{75}\}\subseteq C_1\cap
U_{d_1,2}^F$.
\end{cor}

\begin{proof} Immediate from Proposition~\ref{prop1}, Remark~\ref{rem1}
and Lemma~\ref{lem4}.
\end{proof}

We now need some explicit information concerning the restriction
of characters from $W$ to $W_1$. Using the {\sf CHEVIE} function
{\tt InductionTable}, we compute that 
\begin{align*}
\Res_{W_1}^{W}(512_a)\otimes \varepsilon&=([21,3]\boxtimes {\bf 1})+ 
\mbox{sum of $\psi$ where $\psi\in \Irr(W_1)$ and $a_\psi>4$},\\
\Res_{W_1}^{W}(512_a')\otimes \varepsilon &=([2,31]\boxtimes {\bf 1})+
\mbox{sum of $\psi$ where $\psi\in \Irr(W_1)$ and $a_\psi>4$}.
\end{align*}
Here, ${\bf 1}$ denotes the unit character on the $A_1$-factor of $W_1$ and 
$\varepsilon$ denotes the sign character of $W_1$. The characters of the 
$D_6$-factor are denoted by $[\lambda,\mu]$ where $\lambda$ and $\mu$ 
are partitions such that $|\lambda|+|\mu|=6$. The $a$-invariant 
of a character is defined as in Lusztig \cite{Lu4}, (4.1); in {\sf CHEVIE},
these $a$-invariants are obtained by the function 
{\tt LowestPowerGenericDegrees}.  We have
\[ a_\psi=4 \qquad \mbox{for $\psi=[21,3]\boxtimes {\bf 1}$ and
$\psi=[2,31]\boxtimes {\bf 1}$}.\]
With these explicit formulas, we can now prove the following result.

\begin{prop}\label{final} Assume that the characteristic $p$ is large
enough, such that Lusztig's formula in \cite{LuUS}, Theorem~7.5, for the 
values of a generalized Gelfand--Graev holds for $\Gamma_u^1$. By
\cite{mythom}, Corollary~4.3, there exists some $u \in \{y_{74},y_{75}\}$
such that $\langle E_7[\pm \xi],\Gamma_u \rangle_{G^F}=1$. For this
element $u$, we have 
\[ \big\langle E_7[\pm \xi],\tilde{\Gamma}_u-\tilde{\Gamma}_u' 
\big \rangle_{G^F}=-1.\]
\end{prop}

\begin{proof} We have already mentioned in the remarks preceding
Lemma~\ref{lem4} that 
\[ E_7[\xi]+E_7[-\xi]=R_{512_a}-R_{512_a'}.\]
Since $\tilde{\Gamma}_u$ and $\tilde{\Gamma}_u'$ are rational-valued 
(see Proposition~\ref{lem1}), we have 
\begin{align*}
\big\langle E_7[\pm \xi],\tilde{\Gamma}_u-\tilde{\Gamma}_u' 
\big \rangle_{G^F}
&= \frac{1}{2}\big\langle E_7[\xi]+E_7[-\xi],\tilde{\Gamma}_u-
\tilde{\Gamma}_u' \big \rangle_{G^F}\\
&=\frac{1}{2} \big\langle R_{512_a}-R_{512_a'},\tilde{\Gamma}_u-
\tilde{\Gamma}_u' \big \rangle_{G^F}.
\end{align*}
Now let $\psi\in \Irr(W_1)$ be a constituent in the restriction of
$512_a$ or $512_a'$ from $W$ to $W_1$. Then, by Corollary~\ref{cor1}, 
we must compute the scalar  product $\langle R_\psi,\Gamma_u^1
\rangle_{G_1^F}$.
Let $D$ denote the Alvis--Curtis--Kawanaka duality operation on the 
character ring of $G_1^F$; see Lusztig \cite{Lu4}, (6.8). We have
$D(R_\psi)=R_{\psi \otimes \varepsilon}$ and so 
\[ \langle R_\psi,\Gamma_u^1\rangle_{G_1^F}=
\langle D(R_\psi),D(\Gamma_u^1)\rangle_{G_1^F}=
\langle R_{\psi\otimes \varepsilon},D(\Gamma_u^1)\rangle_{G_1^F}.\]
Now, in order to evaluate the above scalar product, it is enough to
know the values of $R_{\psi \otimes \varepsilon}$ on the unipotent
elements of $G_1^F$. By Shoji's algorithm \cite{Sh1} and by \cite{LuUS},
Corollary~10.9, we know that $R_{\psi\otimes \varepsilon}(y)=0$ if 
$\dim \fB_y^1<a_{\psi \otimes \varepsilon}$. On the other hand, we 
have $D(\Gamma_u^1)(y) =0$ if $\dim \fB_u^1<\dim \fB_y^1$. (This follows
from \cite{LuUS}; see the remarks in \cite{myert03}, (2.4).) Thus, the above 
scalar product is zero if $a_{\psi\otimes\varepsilon}>\dim \fB_u^1=4$. 
Taking into account the explicit information concerning the restrictions
of $512_a$ and $512_a'$ from $W$ to $W_1$, we conclude that 
\begin{equation*}
\big\langle E_7[\pm \xi],\tilde{\Gamma}_u-\tilde{\Gamma}_u' \big 
\rangle_{G^F}=\frac{1}{2}\big\langle R_{[21,3] \boxtimes {\bf 1}}-
R_{[2,31]\boxtimes {\bf 1}},D(\Gamma_u^1)\big \rangle_{G_1^F}.\tag{1}
\end{equation*}
Now $[21,3]\boxtimes {\bf 1}$ and $[2,31]\boxtimes {\bf 1}$ lie in 
the same family of characters of $W_1$; see \cite{Lu4}, Chapter~4. The 
Fourier matrix (which has size $4\times 4$) for that family shows that
\[R_{[21,3]\boxtimes {\bf 1}}-R_{[2,31]\boxtimes {\bf 1}}=-\rho_1-\rho_2\]
where $\rho_1$ and $\rho_2$ are unipotent characters of $G_1^F$. Now, we 
can explicitly compute the unipotent support of these two characters; see 
\cite{LuUS}, \S 11, or \cite{gema}, \S 3.C.  This involves the knowledge 
of the Springer correspondence for $G_1$.  Using the description of that 
correspondence  in \cite{Ca2}, \S 13.3, we find that $\rho_1$ and 
$\rho_2$ have unipotent support $C_1$. Thus, by the formula in 
\cite{gema}, Remark~3.8, we have
\begin{equation*}
\big\langle \rho_i,D(\Gamma_{y_{74}}^1)+D(\Gamma_{y_{75}}^1)
\big\rangle_{G_1^F}=\big\langle D(\rho_i),\Gamma_{y_{74}}^1+
\Gamma_{y_{75}}^1\big\rangle_{G_1^F}=1 \quad \mbox{for $i=1,2$}.\tag{2}
\end{equation*}
Note that $C_{G_1}(y_{74})/C_{G_1}(y_{74})^\circ\cong {\Z}/2{\Z}$ by
Lemma~\ref{lem2} and that $D(\rho_1)$, $D(\rho_2)$ are actual characters 
in the present situation; see \cite{Lu4}, (6.8.2). Now we have $u\in \{y_{74},
y_{75}\}$ and we would like to show that 
\begin{equation*}
\big\langle D(\rho_i),\Gamma_{u}^1\big\rangle_{G_1^F}=
\big\langle \rho_i,D(\Gamma_{u}^1)\big\rangle_{G_1^F}=1 \quad \mbox{for
$i=1,2$}.\tag{3}
\end{equation*}
This can be seen as follows. Fix $i\in \{1,2\}$. Since $D(\rho_i)$ is an
actual character, we certainly have $\langle D(\rho_i),\Gamma_u^1
\rangle_{G_1^F} \geq 0$. Hence, using (2), the latter scalar product equals 
$0$ or $1$. Assume, if possible, that the scalar product is zero. Then the 
scalar product of $-\rho_1- \rho_2$ with $D_G(\Gamma_u^1)$ would be $-1$ 
or $0$. Consequently, the scalar product in (1) would be $-1/2$ or $0$. 
Thus, the only possibility is that the scalar product in (1) equals~$0$. 
But this would mean that 
\[ \big\langle E_7[\pm \xi],\tilde{\Gamma}_u+ \tilde{\Gamma}_u'\big
\rangle_{G^F}=[U_{d,1}^F:U_{d,2}^F]^{1/2}\,\big\langle E_7[\pm \xi],
\Gamma_u\big\rangle_{G^F}=[U_{d,1}^F:U_{d,2}^F]^{1/2}\]
is an even number, which is not true.  So, our assumption was wrong and 
(3) holds. Inserting this into (1), we obtain the desired result.
\end{proof}

\section{Proof of Theorem~\ref{thm1}} 
By \cite{mythom}, Corollary~4.3, the Schur index of $E_7[\pm \xi]$ is at 
most~$2$. Hence, we only need to show that $E_7[\pm \xi]$ cannot be realized
over ${\Q}(\xi)$. Now, we have 
\[\langle E_7[\pm\xi],\Gamma_u\rangle_{G^F}=1\qquad \mbox{for suitable 
$u\in \{y_{74},y_{75}\}$}.\]
So, using the formulas in Definition~\ref{def1},  we obtain that 
\[ \langle E_7[\pm \xi],\tilde{\Gamma}_u+\tilde{\Gamma}_u'\rangle_{G^F}=
[U_{d,1}^F:U_{d,2}^F]^{1/2}=q^m \quad \mbox{for some $m\geq 1$}.\]  
Combining this with Proposition~\ref{final} and using Frobenius reciprocity,
this yields
\[\Big\langle \Res_{U_{d,2}^F.H}^{G^F}\bigl(E_7[\pm \xi]\bigr),\psi_u'
\Big\rangle_{U_{d,2}^F.H}=\big\langle E_7[\pm \xi],\tilde{\Gamma}_u'
\big\rangle_{G^F}=\frac{1}{2}(q^m+1).\]
Since $p\equiv 1 \bmod 4$, we also have $q\equiv 1 \bmod 4$ and so the 
above scalar product is an odd number. Now assume, if possible, that 
$E_7[\pm \xi]$ can be realized over $\Q(\xi)$.  Then the restriction of 
$E_7[\pm \xi]$ to $U_{d,2}^F.H$ can also be realized over ${\Q}(\xi)$. Thus, 
by a standard argument on Schur induces (\cite{Isa}, Corollary~10.2), the 
Schur index of $\psi_u'$ over ${\Q}(\xi)$ divides the above odd number. 
Since the Schur index of $\psi_u'$ over ${\Q}(\xi)$ is at most $2$ (see 
Proposition~\ref{lem1}), it must be one. Thus, $\psi_u'$ can be realized 
over ${\Q}(\xi)$. Now, since $q$ is a square, we have $\xi=\sqrt{-1}$. 
Furthermore, since $p\equiv 1 \bmod 4$, we have $\sqrt{-1} \in {\Q}_p$ 
(the field of $p$-adic numbers). Hence $\psi_u'$ can be realized over 
${\Q}_p$, contradicting Proposition~\ref{lem1}(c). Thus, our assumption 
was wrong and so $E_7[\pm\xi]$ cannot be realized over ${\Q}(\xi)$.

\bigskip
\noindent {\bf Acknowledgements.} I wish to thank Gerhard Hiss for a
careful reading of the manuscript.

\end{document}